\documentclass[11pt,a4paper]{article}

\usepackage[T1]{fontenc}
\usepackage[utf8]{inputenc}
\usepackage{lmodern}
\usepackage[margin=1in]{geometry}

\usepackage{amsmath}
\usepackage{amssymb}
\usepackage{amsfonts}
\usepackage{amsthm}
\usepackage{mathtools}

\usepackage{authblk}

\usepackage[numbers,sort&compress]{natbib}
\bibpunct[, ]{[}{]}{,}{n}{,}{,}

\usepackage{xcolor}
\usepackage[
colorlinks=true,
linkcolor=blue,
citecolor=blue,
urlcolor=blue
]{hyperref}

\theoremstyle{plain}

\newtheorem{theorem}{Theorem}[section]
\newtheorem{lemma}{Lemma}[section]
\newtheorem{corollary}{Corollary}[section]
\newtheorem{proposition}{Proposition}[section]

\theoremstyle{definition}

\newtheorem{definition}{Definition}[section]

\newtheorem{example}{Example}[section]

\theoremstyle{remark}

\newtheorem{remark}{Remark}[section]

\newlength{\tblwidth}
\AtBeginDocument{
	\setlength{\abovedisplayskip}{5pt}
	\setlength{\belowdisplayskip}{5pt}
	\setlength{\abovedisplayshortskip}{3pt}
	\setlength{\belowdisplayshortskip}{3pt}
}
\title{\textbf{Solution-Set Properties of Absolute Value Equations: Inverse Nonnegativity, Boundedness, and Connectedness}}

\author[1]{Abhishek Kumar Singh Sengar%
	\thanks{Corresponding author. Email:
		\href{mailto:abhishek47sengar@gmail.com}
		{abhishek47sengar@gmail.com}}}

\author[2]{Bharat Kumar%
	\thanks{Email:
		\href{mailto:bharatnishad.kanpu@gmail.com}
		{bharatnishad.kanpu@gmail.com}}}

\author[1]{Deepmala%
	\thanks{Email:
		\href{mailto:dmrai23@gmail.com}
		{dmrai23@gmail.com}}}

\affil[1]{Discipline of Natural Sciences,
	Indian Institute of Information Technology Design and Manufacturing
	Jabalpur, Madhya Pradesh 482005, India}

\affil[2]{Department of Economics,
	Indian Institute of Technology Kanpur,
	Uttar Pradesh 208016, India}

\date{}

\begin{document}
	% ================================================================
	\maketitle

\begin{abstract}
	\noindent
The solution-set properties of the absolute value equation (AVE)
\(Ax+|x|=b\) were systematically investigated by Hladík
(\emph{SIAM J. Matrix Anal. Appl.} \textbf{44}, 175-195, 2023),
who also formulated several related open questions. This paper settles and extends some of those questions. A complete converse characterization of inverse nonnegativity is established without restrictions on the matrix order
	or symmetry. The characterization is accompanied by equivalent
	criteria involving nonsingular \(M\)-matrices and spectral radii.
	A diagonal semidefinite certificate is then developed for the uniform
	boundedness of solution sets. Unlike regularity conditions, this
	certificate permits infinitely many solutions, as illustrated by a
	three-dimensional example whose solution set contains a nontrivial
	line segment. Connectedness of the solution set for an individual
	right-hand side is characterized through a finite graph formed from
	its orthant components. A further sufficient condition is obtained
	from the \(P_0\)-property of the associated complementarity matrix.
	The inverse-nonnegativity characterization is finally extended to
	generalized absolute value equations (GAVE) having a nonnegative
	absolute-value coefficient matrix.
\end{abstract}



\section{Introduction}\label{sec:introduction}

For \(A\in\mathbb{R}^{n\times n}\), define the piecewise linear mapping
\(F_A:\mathbb{R}^n\rightarrow\mathbb{R}^n\) by
\(F_A(x)=Ax+|x|\), where the absolute value is applied to each component
of \(x\). Given \(b\in\mathbb{R}^n\), the equation
\begin{equation}\label{eq:ave}
	F_A(x)=b,
	\qquad\text{equivalently,}\qquad
	Ax+|x|=b,
\end{equation}
is referred to as an absolute value equation (AVE). Connections between AVEs,
linear complementarity problems, and mathematical programming have motivated
substantial theoretical and computational research; see
\cite{MangasarianMeyer2006,CottlePangStone2009}.

Existence, uniqueness, reformulations, error estimates, and numerical solution
methods have received considerable attention in the AVE literature; see, for
example,
\cite{MangasarianMeyer2006,WuLi2018,ZamaniHladik2021,ZamaniHladik2023}.
Recent surveys and developments concerning unique solvability and computational
methods are available in \cite{KumarEtAl2024,HladikEtAl2026}.

A systematic investigation of the entire solution set was undertaken by Hladík
\cite{Hladik2023}. For a given right-hand side, write
\begin{equation}\label{eq:solution-set}
	\Sigma_A(b)=\{x\in\mathbb{R}^n:Ax+|x|=b\}.
\end{equation}
Questions considered in \cite{Hladik2023} include whether this set is empty,
finite, bounded, convex, or connected, and whether its members are nonnegative.
Each such property gives rise to a corresponding class of coefficient matrices.
Some of these classes can be recognized efficiently, whereas others lead to
computationally difficult recognition problems.

One unresolved question in \cite{Hladik2023} concerns the interval
\([A-I,A+I]\). It was proved in
\cite[Proposition~3.7]{Hladik2023} that inverse nonnegativity of this interval
guarantees that, for every \(b\geq0\), equation \eqref{eq:ave} has exactly one
solution in \(\mathbb{R}^n\), and this solution is nonnegative. Whether the reverse
implication holds was left unanswered. A converse was obtained in
\cite[Proposition~3.10]{Hladik2023} after additionally assuming regularity of
the interval.

Boundedness and connectedness produced two further questions. By
\cite[Proposition~5.1]{Hladik2023}, boundedness of \(\Sigma_A(b)\) for all
right-hand sides is equivalent to the absence of a nonzero solution of the
corresponding homogeneous AVE. This criterion is exact, but it does not by
itself provide a readily verifiable matrix condition that allows bounded
solution sets with infinitely many elements. A matrix characterization ensuring
connectedness of every nonempty solution set was likewise not obtained.

Related questions have subsequently been considered by Yadav and Dubey \cite{YadavDubey2026a,YadavDubey2026b}.
The alternative equation \(Ax-|x|=b\) was examined in
\cite{YadavDubey2026a}, with emphasis on solvability, convexity, nonnegativity,
and associated matrix classes. The original form \eqref{eq:ave} was considered
in \cite{YadavDubey2026b}. In particular, the converse of
\cite[Proposition~3.7]{Hladik2023} was verified there for matrices of order two
and for symmetric matrices. Its validity for unrestricted square matrices
remained to be determined.

The arguments below make use of linear complementarity theory. For
\(M\in\mathbb{R}^{n\times n}\) and \(q\in\mathbb{R}^n\), the linear
complementarity problem associated with the pair \((M,q)\) is to determine
\(z\in\mathbb{R}^n\) such that
\begin{equation}\label{eq:lcp}
	z\geq0,\qquad Mz+q\geq0,\qquad z^T(Mz+q)=0.
\end{equation}
This problem is denoted by \(\operatorname{LCP}(M,q)\).
If \(A+I\) is nonsingular, the usual transformation of \eqref{eq:ave} leads to
an LCP with
\begin{equation}\label{eq:intro-lcp-matrix}
	M=(A+I)^{-1}(A-I),
	\qquad
	q=(A+I)^{-1}b.
\end{equation}
The correspondence induced by this transformation is developed at the level
of solution sets in Section~\ref{sec:ave-lcp-correspondence}. It permits
properties of \(\Sigma_A(b)\) to be studied through complementarity theory.
Background on the AVE-LCP transformation can be found in
\cite{MangasarianMeyer2006,Prokopyev2009,YadavDubey2026b}.

The main results established in this paper are described below.

\begin{enumerate}
\item The open converse associated with \cite[Proposition~3.7]{Hladik2023} is resolved for square matrices of arbitrary order. Equivalent characterizations in terms of inverse nonnegativity, nonsingular \(M\)-matrices, and spectral radius are
established. The results previously obtained in \cite{YadavDubey2026b} for two-dimensional and symmetric matrices are recovered as special cases.
	
	\item A new sufficient condition is derived under which the solution set
	is bounded for every right-hand side. The condition does not impose
	finiteness. A three-dimensional example is presented for which all
	solution sets are bounded and one of them contains a nontrivial line
	segment.
	
	\item Connectedness for a fixed equation is characterized by a finite
	graph constructed from its orthant solution sets. A reduced graph
	involving only adjacent orthants is also obtained, with its edges
	determined through linear feasibility problems. In addition, a
	structural sufficient condition for connectedness is derived from the
	\(P_0\)-property of the associated complementarity matrix.
	
	\item The inverse-nonnegativity characterization is extended to
	generalized absolute value equations with a nonnegative absolute-value
	coefficient matrix. Corresponding nonsingular \(M\)-matrix and
	spectral-radius characterizations are also established.
\end{enumerate}

The paper is arranged as follows. Preliminary material is presented in
Section~\ref{sec:preliminaries}, and the correspondence between the AVE and
LCP solution sets is developed in
Section~\ref{sec:ave-lcp-correspondence}. Inverse nonnegativity is treated in
Section~\ref{sec:inverse-nonnegative}, while boundedness and connectedness are
examined in Sections~\ref{sec:boundedness} and
\ref{sec:connectedness}, respectively. The generalized absolute value equation is considered in Section~\ref{sec:gave-extension}. Section~\ref{sec:conclusion} concludes the paper.

	\section{Preliminaries}\label{sec:preliminaries}
	
	Throughout the paper, vector and matrix inequalities are understood
	componentwise. The identity matrix of order \(n\) is denoted by \(I\), and
	\(e=(1,\ldots,1)^T\). For a matrix \(C\), its spectral radius is denoted by
	\(\rho(C)\). The nullspace of \(C\) is denoted by \(\ker(C)\).
	
	For \(x\in\mathbb{R}^n\), its positive and negative parts are defined
	componentwise by
	\begin{equation*}
		x_i^+=\max\{x_i,0\},
		\qquad
		x_i^-=\max\{-x_i,0\},
		\qquad i=1,\ldots,n.
	\end{equation*}
	Consequently,
	\begin{equation*}
		x=x^+-x^-,
		\qquad
		|x|=x^++x^-,
		\qquad
		(x^+)^Tx^-=0.
	\end{equation*}
	
	For two matrices \(\underline C,\overline C\in\mathbb{R}^{n\times n}\)
	satisfying \(\underline C\leq\overline C\), the interval matrix
	\([\underline C,\overline C]\) is the set
	\begin{equation*}
		[\underline C,\overline C]
		=
		\{C\in\mathbb{R}^{n\times n}:
		\underline C\leq C\leq\overline C\}.
	\end{equation*}
	
	\begin{definition}
		An interval matrix \([\underline C,\overline C]\) is called regular if
		every \(C\in[\underline C,\overline C]\) is nonsingular. It is called
		inverse nonnegative if every \(C\in[\underline C,\overline C]\) is
		nonsingular and \(C^{-1}\geq0\).
	\end{definition}
	
The following endpoint characterization of inverse nonnegative interval
matrices will be used in Sections~\ref{sec:inverse-nonnegative}
and~\ref{sec:gave-extension}.
	
	\begin{lemma}[Kuttler \cite{Kuttler1971}]
		\label{lem:kuttler}
		The interval matrix \([\underline C,\overline C]\) is inverse
		nonnegative if and only if \(\underline C\) and \(\overline C\) are
		nonsingular and
		\begin{equation*}
			\underline C^{-1}\geq0
			\qquad\text{and}\qquad
			\overline C^{-1}\geq0.
		\end{equation*}
	\end{lemma}
	
	We next recall the matrix classes required in the sequel.
	
\begin{definition}
	Let \(C=(c_{ij})\in\mathbb{R}^{n\times n}\).
	\begin{enumerate}
		\item The matrix \(C\) belongs to the class of \(Z\)-matrices when
		\(c_{ij}\leq0\) for every pair \(i\neq j\).
		
		\item The matrix \(C\) is a \(P\)-matrix when each of its principal
		minors is strictly positive. If these minors are instead required
		to be nonnegative, then \(C\) is a \(P_0\)-matrix.
		
		\item A \(Z\)-matrix \(C\) is a nonsingular \(M\)-matrix when there
		exist a nonnegative matrix \(N\) and a real number \(s\) satisfying
		\(s>\rho(N)\) such that \(C=sI-N\).
	\end{enumerate}
\end{definition}
	
	\begin{definition}
		A matrix \(C\in\mathbb{R}^{n\times n}\) is called a \(Q\)-matrix if
		\(\operatorname{LCP}(C,q)\) has a solution for every
		\(q\in\mathbb{R}^n\).
	\end{definition}
	
	The following standard equivalences for \(Z\)-matrices will be used repeatedly;
	see \cite{BermanPlemmons1994,CottlePangStone2009}.
	
	\begin{lemma}\label{lem:z-matrix-equivalences}
		Let \(C\in\mathbb{R}^{n\times n}\) be a \(Z\)-matrix. Then the
		following statements are equivalent:
		\begin{enumerate}
			\item \(C\) is a \(P\)-matrix;
			\item \(C\) is a \(Q\)-matrix;
			\item \(C\) is a nonsingular \(M\)-matrix;
			\item \(C\) is nonsingular and \(C^{-1}\geq0\).
		\end{enumerate}
	\end{lemma}
	
	The set of all solutions of \(\operatorname{LCP}(M,q)\) defined in \eqref{eq:lcp} is denoted by
	\begin{equation} \label{eq:lcp-solution-set}
		\operatorname{SOL}(M,q)
		=
		\left\{
		z\in\mathbb{R}^n:
		z\geq0,\ Mz+q\geq0,\ z^T(Mz+q)=0
		\right\}.
	\end{equation}
	
The next existence result is due to Karamardian
\cite{Karamardian1976}. We use the formulation stated in
\cite[Theorem~1.1]{SivakumarEtAl2021}.
	\begin{lemma}\label{lem:karamardian}
		Let \(M\in\mathbb{R}^{n\times n}\). Suppose that there exists a vector
		\(d>0\) such that
		\begin{equation*}
			\operatorname{SOL}(M,0)=\{0\}
			\qquad\text{and}\qquad
			\operatorname{SOL}(M,d)=\{0\}.
		\end{equation*}
		Then \(M\) is a \(Q\)-matrix.
	\end{lemma}
	
	We will also use the following characterization of boundedness established by
	Hladík.
	
	\begin{lemma}[{\cite[Proposition~5.1]{Hladik2023}}]
		\label{lem:boundedness-characterization}
		The solution set \(\Sigma_A(b)\) is bounded for every
		\(b\in\mathbb{R}^n\) if and only if the homogeneous equation
		\begin{equation}\label{eq:homogeneous-ave}
			Ax+|x|=0
		\end{equation}
		has only the zero solution.
	\end{lemma}
	
	Finally, we recall a connectedness result for linear complementarity problems.
	
\begin{lemma}[{\cite[Corollary~7]{JonesGowda1998}}]
	\label{lem:p0-connectedness}
	Let \(M\in\mathbb{R}^{n\times n}\) be a \(P_0\)-matrix. If
	\(\operatorname{SOL}(M,q)\) is nonempty and bounded for some
	\(q\in\mathbb{R}^n\), then \(\operatorname{SOL}(M,q)\) is connected.
\end{lemma}

\section{A solution-set correspondence with the linear complementarity problem}
\label{sec:ave-lcp-correspondence}

The algebraic equivalence between absolute value equations and linear
complementarity problems is well known; see
\cite{MangasarianMeyer2006,Prokopyev2009,YadavDubey2026b}. In this section, we
record this equivalence at the level of solution sets. This formulation will
allow boundedness, finiteness, and connectedness to be transferred between an
AVE and its associated LCP.

Suppose that \(A+I\) is nonsingular, and set
\begin{equation}\label{eq:R-definition}
	R=(A+I)^{-1}.
\end{equation}
The matrix and vector in \eqref{eq:intro-lcp-matrix} can then be written as
\begin{equation}\label{eq:M-q-via-R}
	M=R(A-I)=I-2R,
	\qquad
	q=Rb.
\end{equation}

\begin{proposition}\label{prop:ave-lcp-homeomorphism}
	Let \(A+I\) be nonsingular, and let \(R\), \(M\), and \(q\) be
	defined by \eqref{eq:R-definition} and \eqref{eq:M-q-via-R}.
	Then the mapping
	\begin{equation*}
		\Phi:\Sigma_A(b)\longrightarrow\operatorname{SOL}(M,q),
		\qquad
		\Phi(x)=x^-,
	\end{equation*}
	where \(\Sigma_A(b)\) and \(\operatorname{SOL}(M,q)\) are defined
	in \eqref{eq:solution-set} and \eqref{eq:lcp-solution-set},
	respectively, is a homeomorphism. Its inverse is given by
	\(\Phi^{-1}(z)=Mz+q-z\).
\end{proposition}
\begin{proof}
	Let \(x\in\Sigma_A(b)\), and put \(w=x^+\) and \(z=x^-\). Then
	\(x=w-z\), \(|x|=w+z\), \(w\geq0\), \(z\geq0\), and \(w^Tz=0\).
	Substitution into \eqref{eq:ave} gives
	\((A+I)w+(I-A)z=b\). Multiplying this equality by
	\(R=(A+I)^{-1}\), we obtain
	\(w+R(I-A)z=Rb\). Since \(R(I-A)=-M\) and \(Rb=q\), it follows
	that
$
		w=Mz+q.
$
	Hence \(z\geq0\), \(Mz+q\geq0\), and
	\(z^T(Mz+q)=0\). Therefore,
	\(z\in\operatorname{SOL}(M,q)\).
	
	Conversely, let \(z\in\operatorname{SOL}(M,q)\), set
	\(w=Mz+q\), and define \(x=w-z\). Then \(w\geq0\), \(z\geq0\),
	and \(w^Tz=0\). It follows componentwise that
	\(x^+=w\), \(x^-=z\), and \(|x|=w+z\). Moreover,
	\eqref{eq:M-q-via-R} and \(w=Mz+q\) give
	\((A+I)w=(A-I)z+b\). Therefore,
	\((A+I)w+(I-A)z=b\), and hence
	\[
	Ax+|x|=A(w-z)+(w+z)=b.
	\]
	Thus \(x\in\Sigma_A(b)\), and
	\(x=w-z=Mz+q-z\).
	
	The two mappings are mutually inverse. The map \(x\mapsto x^-\) is
	continuous, and its inverse \(z\mapsto Mz+q-z\) is affine. Hence
	\(\Phi\) is a homeomorphism.
\end{proof}

\begin{corollary}\label{cor:solution-set-transfer}
	Under the assumptions of
	Proposition~\ref{prop:ave-lcp-homeomorphism}, the sets
	\(\Sigma_A(b)\) and \(\operatorname{SOL}(M,q)\) have the same
	cardinality. Moreover, \(\Sigma_A(b)\) is nonempty, bounded, finite,
	or connected if and only if \(\operatorname{SOL}(M,q)\) has the
	corresponding property.
\end{corollary}

\begin{proof}
	The assertions concerning cardinality, nonemptiness, finiteness, and
	connectedness follow directly from
	Proposition~\ref{prop:ave-lcp-homeomorphism}. For boundedness, the map
	\(x\mapsto x^-\) sends bounded sets to bounded sets, while its affine
	inverse \(z\mapsto Mz+q-z\) has the same property.
\end{proof}

The following consequence will be used in the inverse-nonnegativity analysis.

\begin{corollary}\label{cor:unique-nonnegative-lcp}
	Under the assumptions of
	Proposition~\ref{prop:ave-lcp-homeomorphism}, suppose that \(q\geq0\).
	Then \eqref{eq:ave} has exactly one solution in \(\mathbb{R}^n\), and
	this solution is nonnegative, if and only if
	\(\operatorname{SOL}(M,q)=\{0\}\). In this case, the unique solution
	of \eqref{eq:ave} is \(x=q\).
\end{corollary}

\begin{proof}
	Since \(q\geq0\), the vector \(z=0\) belongs to
	\(\operatorname{SOL}(M,q)\). By
	Proposition~\ref{prop:ave-lcp-homeomorphism}, a solution
	\(x\in\Sigma_A(b)\) is nonnegative if and only if
	\(\Phi(x)=x^-=0\). Therefore, \eqref{eq:ave} has exactly one solution
	and this solution is nonnegative if and only if \(z=0\) is the only
	element of \(\operatorname{SOL}(M,q)\). The inverse mapping in
	Proposition~\ref{prop:ave-lcp-homeomorphism} then gives \(x=q\).
\end{proof}

\section{Inverse nonnegativity and nonnegative solutions}
\label{sec:inverse-nonnegative}

Hladík \cite[Proposition~3.7]{Hladik2023} proved that inverse
nonnegativity of \([A-I,A+I]\) guarantees that, for every \(b\geq0\),
equation \eqref{eq:ave} has exactly one solution in \(\mathbb{R}^n\), and
this solution is nonnegative. The converse was established in
\cite[Proposition~3.10]{Hladik2023} under the additional assumption that
\([A-I,A+I]\) is regular. More recently, Yadav and Dubey
\cite{YadavDubey2026b} proved the converse for matrices of order two and for
symmetric matrices. We now establish the converse for an arbitrary real square
matrix and obtain two equivalent matrix characterizations.

\begin{theorem}\label{thm:inverse-nonnegative-characterization}
	Let \(A\in\mathbb{R}^{n\times n}\). The following statements are
	equivalent:
	\begin{enumerate}
		\item\label{item:unique-nonnegative}
		For every \(b\geq0\), \eqref{eq:ave} has exactly one
		solution in \(\mathbb{R}^n\), and this solution is nonnegative.
		
		\item\label{item:inverse-nonnegative}
		The interval matrix \([A-I,A+I]\) is inverse nonnegative.
		
		\item\label{item:spectral-characterization}
		The matrix \(A+I\) is nonsingular, and
		$
			R=(A+I)^{-1}\geq0,
			~
			\rho(R)<\frac12.
	$
		
		\item\label{item:mmatrix-characterization}
		The matrix \(A+I\) is nonsingular,
		\(R=(A+I)^{-1}\geq0\), and
		$
			M=I-2R
		$
		is a nonsingular \(M\)-matrix.
	\end{enumerate}
\end{theorem}

\begin{proof}
	We first prove that
	\ref{item:unique-nonnegative} implies
	\ref{item:mmatrix-characterization}. For \(j=1,\ldots,n\), let
	\(x^{(j)}\) be the unique solution of \eqref{eq:ave} corresponding to
	\(b=e_j\), where \(e_j\) is the \(j\)th column of \(I\). By
	\ref{item:unique-nonnegative}, \(x^{(j)}\geq0\). Therefore,
	\(|x^{(j)}|=x^{(j)}\), and
	\begin{equation*}
		(A+I)x^{(j)}=e_j.
	\end{equation*}
	Let \(X=(x^{(1)},\ldots,x^{(n)})\). Then \(X\geq0\) and
	\((A+I)X=I\). Hence \(A+I\) is nonsingular and
$
		R=(A+I)^{-1}=X\geq0.
$
	
	Define \(M=I-2R\). Since \(R\geq0\), every off-diagonal entry of
	\(M\) is nonpositive. Thus \(M\) is a \(Z\)-matrix.
	
	For \(b=0\), the unique solution of \eqref{eq:ave} is \(x=0\).
	By Corollary~\ref{cor:unique-nonnegative-lcp},
	\begin{equation}\label{eq:homogeneous-lcp-zero}
		\operatorname{SOL}(M,0)=\{0\}.
	\end{equation}
	Now choose \(b=e\), where \(e=(1,\ldots,1)^T\), and put
	$
		d=Re.
$
	Because \(R\geq0\) is nonsingular, none of its rows is zero. It follows
	that \(d>0\). Once again, condition
	\ref{item:unique-nonnegative} and
	Corollary~\ref{cor:unique-nonnegative-lcp} give
	\begin{equation}\label{eq:positive-lcp-zero}
		\operatorname{SOL}(M,d)=\{0\}.
	\end{equation}
	Lemma~\ref{lem:karamardian}, applied to
	\eqref{eq:homogeneous-lcp-zero} and
	\eqref{eq:positive-lcp-zero}, shows that \(M\) is a \(Q\)-matrix.
	Since \(M\) is also a \(Z\)-matrix,
	Lemma~\ref{lem:z-matrix-equivalences} implies that \(M\) is a
	nonsingular \(M\)-matrix. Thus
	\ref{item:mmatrix-characterization} holds.
	
	We next show the equivalence of
	\ref{item:mmatrix-characterization} and
	\ref{item:spectral-characterization}. Since \(R\geq0\), the matrix
	\(M=I-2R\) is a nonsingular \(M\)-matrix if and only if
	\begin{equation*}
		1>\rho(2R)=2\rho(R).
	\end{equation*}
	This is equivalent to \(\rho(R)<\frac12\).
	
	We now prove that
	\ref{item:mmatrix-characterization} implies
	\ref{item:inverse-nonnegative}. By
	Lemma~\ref{lem:z-matrix-equivalences},
	\begin{equation*}
		M^{-1}\geq0.
	\end{equation*}
	Moreover,
$
		A-I=(A+I)M.
$
	Consequently,
$
		(A-I)^{-1}=M^{-1}(A+I)^{-1}=M^{-1}R\geq0.
$
We also have
$
		(A+I)^{-1}=R\geq0.
$ 
Lemma~\ref{lem:kuttler} now implies that
	\([A-I,A+I]\) is inverse nonnegative. Hence
	\ref{item:inverse-nonnegative} holds.
	
Finally, suppose that \ref{item:inverse-nonnegative} holds. Since
\([A-I,A+I]\) is inverse nonnegative, it is regular. By
\cite[Theorem~3.3]{WuLi2018}, \eqref{eq:ave} has exactly one
solution for every \(b\in\mathbb{R}^n\). If \(b\geq0\), then
\((A+I)^{-1}\geq0\), and therefore
\(x=(A+I)^{-1}b\geq0\). Since \(|x|=x\), this vector satisfies
\((A+I)x=b\), and thus it is the unique solution of \eqref{eq:ave}.
Hence \ref{item:unique-nonnegative} holds.
\end{proof}

\begin{corollary}\label{cor:explicit-nonnegative-solution}
	Suppose that any of the equivalent conditions in
	Theorem~\ref{thm:inverse-nonnegative-characterization} holds. Then,
	for every \(b\geq0\), the unique solution of \eqref{eq:ave} is
	$
		x=(A+I)^{-1}b.
	$
\end{corollary}

\begin{proof}
	The unique solution is nonnegative. Hence \(|x|=x\), and
	\eqref{eq:ave} reduces to \((A+I)x=b\).
\end{proof}

\begin{remark}\label{rem:uniqueness-interpretation}
	The phrase ``exactly one solution in \(\mathbb{R}^n\), and this
	solution is nonnegative'' in
	Theorem~\ref{thm:inverse-nonnegative-characterization} contains two
	essential requirements.
	
	First, global uniqueness for every \(b\geq0\) is insufficient if the
	corresponding solution is not required to be nonnegative. Yadav and
	Dubey \cite[Example~3.1]{YadavDubey2026b} considered
	\begin{equation*}
		A=
		\begin{pmatrix}
			2&1\\
			1&2
		\end{pmatrix}
	\end{equation*}
	and showed that \(Ax+|x|=b\) has exactly one solution for every
	\(b\geq0\). However, these solutions are not always nonnegative. For
	example, when \(b=(0,1)^T\), the unique solution is
	\(x=(-\frac12,\frac12)^T\). Moreover, for
	\(b=(-1,-1)^T\), every vector \(x\leq0\) satisfying
	\(x_1+x_2=-1\) is a solution. Thus uniqueness for nonnegative
	right-hand sides alone does not imply global unique solvability.
	
	Second, uniqueness only among nonnegative solutions is also
	insufficient. Let \(n=1\) and \(A=0\). For every \(b\geq0\), the
	equation \(|x|=b\) has exactly one nonnegative solution, namely
	\(x=b\). When \(b>0\), however, \(x=-b\) is another solution.
	Furthermore, \([-1,1]\) is not regular and hence is not inverse
	nonnegative.
\end{remark}

\begin{remark}\label{rem:previous-special-cases}
	Theorem~\ref{thm:inverse-nonnegative-characterization} requires neither
	a restriction on the dimension nor symmetry of \(A\). Therefore, the
	two-dimensional and symmetric cases established in
	\cite{YadavDubey2026b} are immediate consequences of the theorem.
\end{remark}

\section{Boundedness without finiteness}
\label{sec:boundedness}

Hladík \cite[Proposition~5.1]{Hladik2023} characterized boundedness of
\(\Sigma_A(b)\) for every right-hand side by the absence of nonzero solutions
to the homogeneous equation \eqref{eq:homogeneous-ave}. The next result provides a practical sufficient condition for this property. The condition permits
\(A-I\) to be singular and therefore does not, in general, imply finiteness of
the solution set.

For a positive diagonal matrix
\(W=\operatorname{diag}(w_1,\ldots,w_n)\), define
\begin{equation}\label{eq:H-W-definition}
	H_W=\frac{WA+A^TW}{2}-W.
\end{equation}
The matrix \(H_W\) is symmetric. The notation \(H_W\succeq0\) means that
\(H_W\) is positive semidefinite, equivalently,
\(y^TH_Wy\geq0\) for every \(y\in\mathbb{R}^n\).

\begin{theorem}\label{thm:boundedness-semidefinite}
	Suppose that there exists a positive diagonal matrix \(W\) such that
	\(H_W\succeq0\) and
	\begin{equation}\label{eq:nullspace-intersection}
		\ker(H_W)\cap\ker(A-I)\cap\mathbb{R}_{-}^{n}
		=\{0\},
	\end{equation}
	where
	\(\mathbb{R}_{-}^{n}=\{x\in\mathbb{R}^n:x\leq0\}\).
	Then \(\Sigma_A(b)\) is bounded for every
	\(b\in\mathbb{R}^n\).
\end{theorem}

\begin{proof}
	By Lemma~\ref{lem:boundedness-characterization}, it is sufficient to
	show that the homogeneous equation
	\begin{equation*}
		Ax+|x|=0
	\end{equation*}
	has only the zero solution.
	
	Let \(x\) be a solution of this equation. Multiplying it from the left
	by \(x^TW\), we obtain
	\begin{equation*}
		0=x^TWAx+x^TW|x|.
	\end{equation*}
	Since \(x^TWAx\) is a scalar,
	\begin{equation*}
		x^TWAx
		=
		x^T\left(\frac{WA+A^TW}{2}\right)x.
	\end{equation*}
	Using \eqref{eq:H-W-definition}, it follows that
	\begin{equation}\label{eq:boundedness-energy}
		0
		=
		x^TH_Wx+
		\sum_{i=1}^{n}w_i\bigl(x_i^2+x_i|x_i|\bigr).
	\end{equation}
	
	The first term on the right-hand side of
	\eqref{eq:boundedness-energy} is nonnegative because
	\(H_W\succeq0\). For each \(i\),
	\begin{equation*}
		x_i^2+x_i|x_i|
		=
		\begin{cases}
			2x_i^2, & x_i>0,\\
			0,      & x_i\leq0.
		\end{cases}
	\end{equation*}
	Hence every term on the right-hand side of
	\eqref{eq:boundedness-energy} is nonnegative. Equality is possible
	only if \(x\leq0\) and \(x^TH_Wx=0\). Since \(H_W\) is positive
	semidefinite, \(x^TH_Wx=0\) implies \(H_Wx=0\).
	
	Because \(x\leq0\), we have \(|x|=-x\). The homogeneous equation then
	reduces to
	\begin{equation*}
		(A-I)x=0.
	\end{equation*}
	Thus
	\begin{equation*}
		x\in
		\ker(H_W)\cap\ker(A-I)\cap\mathbb{R}_{-}^{n}.
	\end{equation*}
	Condition \eqref{eq:nullspace-intersection} yields \(x=0\). The result
	now follows from
	Lemma~\ref{lem:boundedness-characterization}.
\end{proof}

A useful specialization is obtained when \(A-I\) is skew-symmetric.

\begin{corollary}\label{cor:skew-boundedness}
	Let \(A=I+K\), where \(K^T=-K\). If
	\begin{equation}\label{eq:skew-kernel-condition}
		\ker(K)\cap\mathbb{R}_{-}^{n}=\{0\},
	\end{equation}
	then \(\Sigma_A(b)\) is bounded for every
	\(b\in\mathbb{R}^n\).
\end{corollary}

\begin{proof}
	Take \(W=I\). Since \(K^T=-K\),
	\begin{equation*}
		H_I
		=
		\frac{A+A^T}{2}-I
		=
		0.
	\end{equation*}
	Moreover, \(A-I=K\). Therefore,
	condition \eqref{eq:nullspace-intersection} reduces to
	\eqref{eq:skew-kernel-condition}, and the conclusion follows from
	Theorem~\ref{thm:boundedness-semidefinite}.
\end{proof}

The following example shows that
Theorem~\ref{thm:boundedness-semidefinite} does not require the solution sets
to be finite.

\begin{example}\label{ex:bounded-infinite}
	Consider
	\begin{equation*}
		K=
		\begin{pmatrix}
			0 & -1 & -1\\
			1 &  0 & -1\\
			1 &  1 &  0
		\end{pmatrix},
		\qquad
		A=I+K=
		\begin{pmatrix}
			1 & -1 & -1\\
			1 &  1 & -1\\
			1 &  1 &  1
		\end{pmatrix}.
	\end{equation*}
	Clearly, \(K^T=-K\), and a direct calculation gives
	\(\ker(K)=\operatorname{span}\{(1,-1,1)^T\}\). This nullspace
	contains no nonzero componentwise nonpositive vector. Hence
	\(\ker(K)\cap\mathbb{R}_{-}^{3}=\{0\}\). By
	Corollary~\ref{cor:skew-boundedness}, \(\Sigma_A(b)\) is bounded for
	every \(b\in\mathbb{R}^3\).
	
	Now take
	\begin{equation*}
		b=
		\begin{pmatrix}
			2\\
			0\\
			-2
		\end{pmatrix},
		\qquad
		x(t)=
		\begin{pmatrix}
			-1\\
			-1\\
			-1
		\end{pmatrix}
		+
		t
		\begin{pmatrix}
			1\\
			-1\\
			1
		\end{pmatrix},
		\qquad -1\leq t\leq1.
	\end{equation*}
	For every \(t\in[-1,1]\), we have \(x(t)\leq0\), and therefore
	\(|x(t)|=-x(t)\). Moreover,
	\(K(1,-1,1)^T=0\) and
	\(K(-1,-1,-1)^T=(2,0,-2)^T=b\). Consequently,
	\begin{equation*}
		Ax(t)+|x(t)|
		=(I+K)x(t)-x(t)
		=Kx(t)
		=b.
	\end{equation*}
	Thus \(x(t)\in\Sigma_A(b)\) for every \(t\in[-1,1]\).
	Therefore, \(\Sigma_A(b)\) is bounded but contains infinitely many
	points.
\end{example}

\begin{remark}\label{rem:boundedness-not-finiteness}
	In Example~\ref{ex:bounded-infinite}, the matrix
	\(A-I=K\) is singular. Hence the example lies outside conditions that
	require all selection matrices associated with the orthants to be
	nonsingular. It demonstrates explicitly that boundedness of every
	solution set does not imply finiteness.
\end{remark}

	\section{Connectedness of the solution set}
	\label{sec:connectedness}
	
	We first describe connectedness by using the orthant decomposition of
	\(\mathbb{R}^n\). For a sign vector
	\(s=(s_1,\ldots,s_n)^T\in\{-1,1\}^n\), let
	\(D_s=\operatorname{diag}(s_1,\ldots,s_n)\), and define
	$$
		P_s(b)
		=
		\{x\in\mathbb{R}^n:(A+D_s)x=b,\ D_sx\geq0\}.
$$
	Each nonempty \(P_s(b)\) is a closed convex polyhedron and is therefore
	connected. Since \(|x|=D_sx\) whenever \(D_sx\geq0\), the complete solution
	set satisfies
	\begin{equation}\label{eq:orthant-union}
		\Sigma_A(b)=\bigcup_{s\in\{-1,1\}^n}P_s(b).
	\end{equation}
	
	Let
	\(\mathcal{S}_A(b)=\{s\in\{-1,1\}^n:P_s(b)\neq\varnothing\}\).
	For a fixed \(A\) and \(b\), define the orthant-intersection graph
	\(\mathcal{G}_A(b)\) as follows. Its vertex set is
	\(\mathcal{S}_A(b)\), and two distinct vertices \(s\) and \(t\) are joined
	by an edge if
	\begin{equation*}
		P_s(b)\cap P_t(b)\neq\varnothing.
	\end{equation*}
	
	\begin{theorem}\label{thm:orthant-graph-connectedness}
		Suppose that \(\Sigma_A(b)\neq\varnothing\). Then
		\(\Sigma_A(b)\) is connected if and only if the graph
		\(\mathcal{G}_A(b)\) is connected.
	\end{theorem}
	
	\begin{proof}
		Suppose first that \(\mathcal{G}_A(b)\) is connected. Each vertex of
		the graph represents a connected set \(P_s(b)\). Along every edge,
		the corresponding sets have a nonempty intersection. It follows by
		successively joining intersecting connected sets that their union is
		connected. By \eqref{eq:orthant-union}, this union is
		\(\Sigma_A(b)\).
		
		Conversely, suppose that \(\mathcal{G}_A(b)\) is disconnected. Its
		vertex set can then be partitioned into two nonempty subsets
		\(\mathcal{S}_1\) and \(\mathcal{S}_2\) with no edge joining them.
		Define
		\begin{equation*}
			\mathcal{U}_k
			=
			\bigcup_{s\in\mathcal{S}_k}P_s(b),
			\qquad k=1,2.
		\end{equation*}
		Since there is no edge between the two vertex sets,
		\(\mathcal{U}_1\cap\mathcal{U}_2=\varnothing\). Each
		\(\mathcal{U}_k\) is a finite union of closed sets and is therefore
		closed. Moreover,
		\(\Sigma_A(b)=\mathcal{U}_1\cup\mathcal{U}_2\), and each
		\(\mathcal{U}_k\) is the complement of the other relative to
		\(\Sigma_A(b)\). Thus both sets are relatively open as well as
		relatively closed in \(\Sigma_A(b)\). Hence \(\Sigma_A(b)\) is
		disconnected.
	\end{proof}
	
	The graph in Theorem~\ref{thm:orthant-graph-connectedness} can be reduced to
	edges between adjacent sign vectors. Two sign vectors are called adjacent if
	they differ in exactly one coordinate.
	
	\begin{corollary}\label{cor:adjacent-orthant-graph}
		Suppose that \(\Sigma_A(b)\neq\varnothing\). Form a graph with vertex
		set \(\mathcal{S}_A(b)\) by joining two vertices \(s\) and \(t\) only
		when they are adjacent and \(P_s(b)\cap P_t(b)\neq\varnothing\).
		Then \(\Sigma_A(b)\) is connected if and only if this reduced graph
		is connected.
	\end{corollary}
	
	\begin{proof}
		It is sufficient to show that every edge of
		\(\mathcal{G}_A(b)\) can be replaced by a path consisting of edges
		between adjacent sign vectors.
		
		Let \(s,t\in\mathcal{S}_A(b)\), and suppose that
		\(x\in P_s(b)\cap P_t(b)\). If \(s_i\neq t_i\), then the inequalities
		\(s_ix_i\geq0\) and \(t_ix_i\geq0\) imply \(x_i=0\). Starting from
		\(s\), change the differing signs one at a time until \(t\) is
		reached. For every intermediate sign vector \(r\), the changed
		coordinates of \(x\) are zero. Consequently,
		\(D_rx=D_sx=|x|\), and hence
		\begin{equation*}
			(A+D_r)x=Ax+|x|=b,
			\qquad
			D_rx\geq0.
		\end{equation*}
		Thus \(x\in P_r(b)\) for every intermediate sign vector. Consecutive
		vertices on the resulting path are adjacent and their corresponding
		polyhedra contain the common point \(x\).
	\end{proof}
	
	\begin{remark}\label{rem:edge-feasibility}
		If adjacent sign vectors \(s\) and \(t\) differ in coordinate \(i\),
		then \(P_s(b)\cap P_t(b)\neq\varnothing\) if and only if the linear
		system
		\begin{equation*}
			(A+D_s)x=b,
			\qquad
			D_sx\geq0,
			\qquad
			x_i=0
		\end{equation*}
		is feasible. Thus the vertices and edges of the reduced graph can be
		determined by finitely many linear feasibility problems. The
		characterization remains exponential in the dimension because as many
		as \(2^n\) orthants may have to be examined.
	\end{remark}
	
	The orthant graph gives an exact connectedness characterization for a fixed
	equation.
	We next derive a structural sufficient condition that applies simultaneously
	to all bounded solution sets.
	
\begin{theorem}\label{thm:p0-connectedness}
	Suppose that \(A+I\) is nonsingular and
	\(M=(A+I)^{-1}(A-I)\) is a \(P_0\)-matrix. If
	\(\Sigma_A(b)\) is nonempty and bounded, then
	\(\Sigma_A(b)\) is connected.
\end{theorem}
	
	\begin{proof}
		Let \(q=(A+I)^{-1}b\). By
		Corollary~\ref{cor:solution-set-transfer},
		\(\operatorname{SOL}(M,q)\) is nonempty and bounded. Since \(M\) is
		a \(P_0\)-matrix, Lemma~\ref{lem:p0-connectedness} shows that
		\(\operatorname{SOL}(M,q)\) is connected. The homeomorphism in
		Proposition~\ref{prop:ave-lcp-homeomorphism} then implies that
		\(\Sigma_A(b)\) is connected.
	\end{proof}
	
	Combining Theorems~\ref{thm:boundedness-semidefinite} and
	\ref{thm:p0-connectedness} yields the following condition for all
	right-hand sides.
	
	\begin{corollary}\label{cor:bounded-connected}
		Suppose that \(A+I\) is nonsingular and
		\(M=(A+I)^{-1}(A-I)\) is a \(P_0\)-matrix. If there exists a positive
		diagonal matrix \(W\) such that \(H_W\succeq0\) and
		\begin{equation*}
			\ker(H_W)\cap\ker(A-I)\cap\mathbb{R}_{-}^{n}
			=\{0\},
		\end{equation*}
		then, for every \(b\in\mathbb{R}^n\), the set \(\Sigma_A(b)\) is
		either empty or bounded and connected.
	\end{corollary}
	
	\begin{proof}
		Theorem~\ref{thm:boundedness-semidefinite} shows that
		\(\Sigma_A(b)\) is bounded for every \(b\). Whenever it is nonempty,
		Theorem~\ref{thm:p0-connectedness} shows that it is connected.
	\end{proof}
	
	The matrix from Example~\ref{ex:bounded-infinite} also illustrates the preceding
	connectedness result.
	
	\begin{example}\label{ex:bounded-connected-infinite}
		Let \(A\) be the matrix from Example~\ref{ex:bounded-infinite}. Then
		\(A+I\) is nonsingular, and direct calculation gives
		\begin{equation*}
			M=(A+I)^{-1}(A-I)
			=
			\frac17
			\begin{pmatrix}
				2 & -1 & -3\\
				3 &  2 & -1\\
				1 &  3 &  2
			\end{pmatrix}.
		\end{equation*}
		The principal minors of order one are all equal to \(2/7\), while
		the principal minors of order two are all equal to \(1/7\). Moreover,
		\(\det(M)=0\). Therefore, \(M\) is a \(P_0\)-matrix but not a
		\(P\)-matrix.
		
		Example~\ref{ex:bounded-infinite} shows that \(\Sigma_A(b)\) is
		bounded for every \(b\in\mathbb{R}^3\). Hence
		Theorem~\ref{thm:p0-connectedness} implies that every nonempty
		\(\Sigma_A(b)\) is connected. For
		\(b=(2,0,-2)^T\), the solution set contains the line segment displayed
		in Example~\ref{ex:bounded-infinite}. Thus this solution set is
		bounded, connected, and infinite.
	\end{example}

\section{Extension to generalized absolute value equations}
\label{sec:gave-extension}

For \(A,B\in\mathbb{R}^{n\times n}\), define the piecewise linear mapping
\(F_{A,B}:\mathbb{R}^n\rightarrow\mathbb{R}^n\) by
\(F_{A,B}(x)=Ax+B|x|\). Given \(b\in\mathbb{R}^n\), the equation
\begin{equation}\label{eq:gave-plus}
	F_{A,B}(x)=b,
	\qquad\text{equivalently,}\qquad
	Ax+B|x|=b,
\end{equation}
is referred to as a generalized absolute value equation (GAVE). The
inverse-nonnegativity characterization obtained for the AVE is now extended
to \eqref{eq:gave-plus}. Throughout this section, \(B\geq0\) is assumed;
consequently, \(A-B\leq A+B\).

The following standard characterization will be used. It is stated here in
the plus-sign form corresponding to \eqref{eq:gave-plus}.

\begin{lemma}[{\cite[Theorem~3.2]{WuShen2021}}]
	\label{lem:gave-plus-criterion}
	The GAVE \eqref{eq:gave-plus} has a unique solution for every
	\(b\in\mathbb{R}^n\) if and only if \(A+BD\) is nonsingular for every
	diagonal matrix \(D=\operatorname{diag}(d_1,\ldots,d_n)\) satisfying
	\(|d_i|\leq1\), \(i=1,\ldots,n\).
\end{lemma}

\begin{theorem}\label{thm:gave-inverse-nonnegative}
	Let \(B\geq0\). The following statements are equivalent:
	\begin{enumerate}
		\item\label{item:gave-unique-nonnegative}
		For every \(b\geq0\), the GAVE \eqref{eq:gave-plus} has exactly one
		solution in \(\mathbb{R}^n\), and this solution is nonnegative.
		
		\item\label{item:gave-inverse-nonnegative}
		The interval matrix \([A-B,A+B]\) is inverse nonnegative.
		
		\item\label{item:gave-spectral}
		The matrix \(A+B\) is nonsingular,
		\(R=(A+B)^{-1}\geq0\), and \(\rho(RB)<\frac12\).
		
		\item\label{item:gave-mmatrix}
		The matrix \(A+B\) is nonsingular,
		\(R=(A+B)^{-1}\geq0\), and
		\(M=R(A-B)=I-2RB\) is a nonsingular \(M\)-matrix.
	\end{enumerate}
\end{theorem}

\begin{proof}
	Suppose first that \ref{item:gave-unique-nonnegative} holds. For each
	\(j=1,\ldots,n\), let \(x^{(j)}\) denote the unique solution of
	\eqref{eq:gave-plus} corresponding to \(b=e_j\). Since
	\(x^{(j)}\geq0\), one has \((A+B)x^{(j)}=e_j\). Define
	\(X=(x^{(1)},\ldots,x^{(n)})\). Then \((A+B)X=I\), and hence \(A+B\)
	is nonsingular and \(X=(A+B)^{-1}\geq0\). Set
	\(R=(A+B)^{-1}\).
	
	Write \(x=u-v\), where \(u=x^+\) and \(v=x^-\). Then \(u\geq0\),
	\(v\geq0\), and \(u^Tv=0\). Substitution into \eqref{eq:gave-plus}
	gives \((A+B)u+(B-A)v=b\). Multiplication by \(R\) yields
	\(u=Mv+Rb\), where \(M=R(A-B)=I-2RB\). Consequently, the solutions
	of \eqref{eq:gave-plus} correspond bijectively to the solutions of
	\(\operatorname{LCP}(M,Rb)\).
	
	Since \(R\geq0\) and \(B\geq0\), the matrix \(M=I-2RB\) is a
	\(Z\)-matrix. For \(b=0\), the unique solution of the GAVE is \(x=0\);
	therefore, \(\operatorname{SOL}(M,0)=\{0\}\). Let
	\(e=(1,\ldots,1)^T\) and \(d=Re\). Since \(R\geq0\) is nonsingular,
	no row of \(R\) is zero, and hence \(d>0\). For \(b=e\), the unique
	GAVE solution is nonnegative, so its negative part is zero. The
	bijection therefore gives \(\operatorname{SOL}(M,d)=\{0\}\).
	Lemma~\ref{lem:karamardian} now implies that \(M\) is a \(Q\)-matrix.
	Since \(M\) is also a \(Z\)-matrix,
	Lemma~\ref{lem:z-matrix-equivalences} shows that \(M\) is a
	nonsingular \(M\)-matrix. Thus \ref{item:gave-mmatrix} holds.
	
	Statements \ref{item:gave-mmatrix} and \ref{item:gave-spectral} are
	equivalent. Indeed, \(RB\geq0\), and \(I-2RB\) is a nonsingular
	\(M\)-matrix if and only if \(\rho(2RB)<1\), which is equivalent to
	\(\rho(RB)<\frac12\).
	
	Next, suppose that \ref{item:gave-mmatrix} holds. Since \(M\) is a
	nonsingular \(M\)-matrix, \(M^{-1}\geq0\). Moreover,
	\(A-B=(A+B)M\), and therefore
	\((A-B)^{-1}=M^{-1}(A+B)^{-1}=M^{-1}R\geq0\). Together with
	\((A+B)^{-1}=R\geq0\), Lemma~\ref{lem:kuttler} implies that
	\([A-B,A+B]\) is inverse nonnegative. Thus
	\ref{item:gave-inverse-nonnegative} holds.
	
	Finally, suppose that \ref{item:gave-inverse-nonnegative} holds. Let
	\(D=\operatorname{diag}(d_1,\ldots,d_n)\), where \(|d_i|\leq1\).
	Since \(B\geq0\), one has \(A-B\leq A+BD\leq A+B\). Hence
	\(A+BD\in[A-B,A+B]\), and inverse nonnegativity of the interval
	implies that \(A+BD\) is nonsingular. By
	Lemma~\ref{lem:gave-plus-criterion}, \eqref{eq:gave-plus} has exactly
	one solution for every \(b\in\mathbb{R}^n\).
	
	Now let \(b\geq0\). Since \(A+B\in[A-B,A+B]\), inverse
	nonnegativity gives \((A+B)^{-1}\geq0\). Thus
	\(x=(A+B)^{-1}b\geq0\), and consequently
	\(Ax+B|x|=(A+B)x=b\). The global uniqueness established above shows
	that \(x\) is the only solution. Therefore,
	\ref{item:gave-unique-nonnegative} holds.
\end{proof}

\begin{corollary}\label{cor:ave-from-gave}
	For \(B=I\), Theorem~\ref{thm:gave-inverse-nonnegative} reduces to
	Theorem~\ref{thm:inverse-nonnegative-characterization}.
\end{corollary}

	\section{Conclusion}\label{sec:conclusion}
The open converse concerning inverse nonnegativity of \([A-I,A+I]\) has been completely resolved for matrices of arbitrary order, thereby removing the dimensional and symmetry restrictions imposed in earlier results. Equivalent descriptions through a nonsingular \(M\)-matrix and a spectral-radius inequality have also been obtained, and the resulting characterization has been extended to generalized absolute value equations with a nonnegative absolute-value coefficient matrix. The boundedness question has been advanced by a diagonal semidefinite condition that accommodates
infinite solution sets, as demonstrated by an explicit three-dimensional
example. For connectedness, an exact graph characterization has been developed
for a fixed equation, together with a \(P_0\)-matrix condition ensuring that a
nonempty bounded solution set is connected. Thus the questions raised in
\cite{Hladik2023} have been resolved or strengthened in several directions,
although a complete matrix characterization guaranteeing connectedness for
every nonempty solution set remains open.
	
\section*{Acknowledgements}
The first author acknowledges the support of the Ministry of Education, Government of India, through the GATE fellowship (Registration No. MA22S25023148).	
	
\subsection*{Data Availability}
No data were used for the research presented in this work.
	
\section*{Statements and Declarations}

\noindent
\textbf{Funding.}
The authors received no funding for this work.

\medskip
\noindent
\textbf{Competing Interests.}
The authors declare that they have no financial or non-financial interests
that are directly or indirectly related to the work reported in this article.

\end{document}